\documentclass[11pt]{amsart}
\usepackage{latexsym}
\usepackage{amsmath,amsthm,amssymb}

\newtheorem{theorem}{Theorem}[section]

\newtheorem{lemma}[theorem]{Lemma}

\newcommand \A {{\mathcal A}}
\newcommand \B {{\mathcal B}}

\newcommand \Z {{\mathbb Z}}

\begin{document}
\title{Quasi-Isometries between groups with infinitely many ends}
\author{Panos Papasoglu, Kevin Whyte}
\date{}

\begin{abstract}
Let $G,F$ be finitely generated groups with infinitely many ends and let 
$\pi_1(\Gamma,\A )$, $\pi_1(\Delta ,\B )$ be graph of groups decompositions
of $F,G$ such that all edge groups are finite and all vertex groups
have at most one end. We show that $G,F$ are quasi-isometric if and
only if every one-ended vertex group of $\pi_1(\Gamma,\A )$ is quasi-isometric
to some one-ended vertex group of  $\pi_1(\Delta ,\B )$ and 
every one-ended vertex group of $\pi_1(\Delta ,\B )$ is quasi-isometric
to some one-ended vertex group of  $\pi_1(\Gamma,\A )$. From our proof it also 
follows that if $G$ is any finitely generated group, of order at least 
three, the groups: $G\ast G, G\ast \Z,G\ast G \ast G$ and 
$G\ast \Z/2\Z$ are all quasi-isometric.
\end{abstract}

\maketitle
\section*{Introduction}
  One of the most appealing and influential theorems in geometric 
group theory is Stallings' Ends Theorem, \cite{St1}, \cite{St2}.  This theorem 
says that a finitely generated group splits as a free product or
HNN-extension with finite  amalgamation if and only if it has more than
one end.  The property of having  infinitely many ends is geometric, in
particular, it is invariant under  quasi-isometry.  One of the main goals
of geometric group theory is to classify finitely generated groups up to
quasi-isometry. Given that splitting over finite subgroups is invariant
under  quasi-isometry, it is natural to ask how the quasi-isometry type
of a free product with finite amalgamation is related to the
types of its factors.

   This question is not as straightforward as one might think.  It is 
not true that if $G$ and $G'$ are quasi-isometric then $G*H$ and 
$G'*H$ are.  Examples of this can already be seen among finite groups: $
\Z/2\Z * { \Z}/2{\Z}$ and ${\Z}/2{\Z} * {\Z}/3{\Z}$ are not 
quasi-isometric.  Our first theorem shows
that this is essentially the only source of examples. 

\begin{theorem}\label{FreeProducts}
Let $A, B$ and $C$ be nontrivial groups.  
If $A$ and $B$ are quasi-isometric then $A*C$ and $B*C$ are 
quasi-isometric unless $C$, and one of $A$ or $B$, are of order $2$.
\end{theorem}

We note that some special cases of theorem \ref{FreeProducts} have been treated
in \cite {P}, and \cite {W}.

Stallings' theorem gives splittings over finite subgroups, 
not free products.  Our next theorem shows that from the 
quasi-isometric point of view, finite amalgamated products are free.

\begin{theorem}\label{Amalgamation}
 Let $A$ and $B$ be groups, and $F$ a common 
finite proper subgroup.  Unless $F$ is of index 2 in both $A$ and $B$ 
then $A*B$ and $A*_{F}B$ are quasi-isometric.  Likewise,  $A*_{F}$ 
and $A * \Z$ are quasi-isometric.
\end{theorem}

  From these theorems we get a complete classification of the 
quasi-isometry types of graphs of groups with finite edge groups:

\begin{theorem}\label{General}
Let $G,H$ be finitely generated groups with infinitely many ends and let
$\pi _1(\Gamma,\A )$,  $\pi _1(\Delta ,\B )$ be decompositions
of $G,H$ in graphs of groups such that all edge groups are finite.
If $\pi _1(\Gamma,\A )$,  $\pi _1(\Delta ,\B )$ have the same set of 
quasi-isometry types of vertex groups (without multiplicities) 
then $G$ and $H$ are quasi-isometric. 
\end{theorem}

  It is natural to ask for a converse.  As it may be possible to split a 
vertex group, the naive converse cannot be true.  By Stallings' theorem, 
a vertex group splits if and only if it has $2$ or more ends.  The obvious thing 
to do is to split until no more splitting is possible, in other words, 
until all vertex groups are finite or one ended.  A finite graph of 
groups with this property is called {\bf terminal}, and a group which has
a terminal  splitting is called {\bf accessible}.

   The Grushko-Neumann theorem shows that finitely generated, 
torsion-free groups, are accessible.  It is also true that 
finitely presented groups are accessible (\cite{D1}).  While it would 
be natural to think all finitely generated groups are accessible,  
indeed this was a conjecture for quite a while, it is not true 
(\cite{D2}).  It follows easily from the characterization of 
accessibility in \cite {TW} that accessibility is a 
quasi-isometry invariant.   Together with the earlier theorems, this
yields:

\begin{theorem}\label{Accessible}  Let $G$ be an accessible group
and let $\pi _1(\Gamma,\A )$ be a terminal graph of groups decomposition 
of $G$. A group $G'$ is quasi-isometric to $G$ if and only if it is 
also accessible and any terminal decomposition of $G'$ , $\pi _1(\Delta ,\B )$, 
has the same set of quasi-isometry types of one ended factors and the same 
number of ends.
\end{theorem}
 
   This can be viewed as a step in Gromov's program (\cite{G}) to 
classify finitely generated groups up to quasi-isometry.   It 
effectively reduces the classification of accessible groups to the 
classification of one-ended groups.  It would be very interesting to 
have a similar reduction for the classification of non-accessible groups,
perhaps to some sort of quasi-conformal structure on the set of ends
together with the quasi-isometry types of one-ended factors.

  The first author would like to thank Pierre de la Harpe , Thomas
Delzant and Frederic Haglund for conversations related to this work. The 
second author would like to thank Benson Farb, Lee Mosher, and 
Shmuel Weinberger for their suggestions and encouragement.
We would also like to thank the referee for many improvements
and corrections.

\section{Basic Construction}

  All of our spaces are the vertex sets of connected graphs of bounded 
valence.  We give these spaces the path metric of the graph, where 
every edge is considered to have length $1$.  The primary motivating 
examples are the Cayley graphs of finitely generated groups.  
Different finite generating sets give different graphs, but the 
induced metrics on the group are bilipschitz equivalent.

  If $X$ is a graph, a {\bf net} in $X$ is a subset $S$ which is {\bf
coarsly dense}, meaning that there is an $r>0$ so that every $x$ in $X$
is within $r$ of some $s \in S$.  The inclusion map $S \to X$ is a
quasi-isometry when $S$ is given the induced metric.  One can give $S$ a
graph structure by connecting any two vertices within $2r$ by an edge. 
The resulting metric is bilipschitz equivalent to the metric induced
from $X$.   

  The free product, $G*H$, of two finitely presented groups has a nice
geometric model.  Bass-Serre theory gives a tree $T$ with a $G*H$ action,
free on edges, with quotient an edge, and the stabilizers of the
vertices the conjugates of $G$ and $H$.  The model for $G*H$ is
produced by "blowing up" the vertices of the tree to be copies of the
Cayley graphs of $G$ and $H$, so that the $G*H$ action becomes free. 
The resulting space has all its vertices in these vertex subgraphs, and
there is exactly one edge at every vertex connecting to another vertex
space.  See \cite{SW} for more details and generalizations to more
complicated graphs of groups.  

  We need to generalize this and define the {\em free product} of two
spaces, $X$ and $Y$.  Much of the tree of spaces structure of free
products of groups makes sense for arbitrary spaces: one wants a graph
with distinguished subgraphs, each isomorphic to $X$ or $Y$, which
are disjoint and cover all the vertices.  Every edge not in one of
these subgraphs should connect a subgraph isomorphic to $X$ to one
isomorphic to $Y$, and there should be precisely one such edge at every
vertex.  Finally, the pattern of attachments of these subgraphs should 
be a tree.  

  This description is not quite sufficient to uniquely define a
graph.  To construct such a tree of spaces, we start with, say, a 
copy of $X$ and, at every vertex of this $X$ add an edge connecting
to a copy of $Y$.  Immediately we run into ambiguity - connecting to
a copy of $Y$ at what point?  This difficulty does not arise when
building a free product out of Cayley graphs because they have 
transitive isomorphism group, which makes all possible points of 
attachment equivalent.

  To get around this, we give our spaces $X$ and $Y$ distinguished
base points, $x_0$ and $y_0$.  We can then construct a canonical 
free product of $(X,x_0)$ and $(Y,y_0)$ as follows:

  Let $\Gamma_{0}$ be the graph which is the disjoint union of $X$ and 
$Y$, with an edge added connecting the base points.  Observe that this 
graph satisfies all the above conditions, except that some vertices 
are not incident to edges not in a subgraph, although the base point 
of every copy of $X$ or $Y$ is.

  Given $\Gamma_{n}$ build $\Gamma_{n+1}$ as follows: For any $v$ in 
$\Gamma_{n}$ not incident to an edge which connects to another 
subgraph, add a new subgraph, isomorphic to $X$ or $Y$ as required, 
and an edge which connects to $v$ and to the base point of the added 
graph.  One has $\Gamma_{n}$ embedded canonically as a subgraph of 
$\Gamma_{n+1}$, and the direct limit (union) as $n \to \infty$ is 
$(X,x_0)*(Y,y_0)$.

  The space $(X,x_0)*(Y,y_0)$ is characterized by:

\begin{itemize}
\item $X*Y$ contains a disjoint collection of subgraphs, each with an 
isomorphism to $X$ or $Y$.
\item Every vertex of $X*Y$ is contained in one of the subgraphs and 
is incident to exactly one edge not in that subgraph.
\item Every edge not in one of the subgraphs connects a subgraph 
isomorphic to $X$ and a subgraph isomorphic to $Y$, and is incident to 
the base point of one of the components it connects.  Further, there 
is a unique edge in $X*Y$, called the base edge, which is incident to 
the base points in both components it connects.
\item The quotient graph in which each of the subgraphs is collapsed 
to a point is a tree.
\end{itemize}

  As discussed above, the construction does not depend on base points
for groups.  For more general graphs as $X$ and $Y$ the choice of base
points will affect the graph constructed by the above.  The bilipschitz
class of metric space is independent of these choices for a wider class
of spaces.

  We say that $X$ is {\bf homogeneous } if $X$ has the property that for
some $L$ and for any $x_{1}$  and $x_{2}$ in $X$ there is a self
$L$-bilipschitz  map taking $x_{1}$ to $x_{2}$.  Note that
this is much weaker than transitive isometry group. 

\begin{lemma}\label{invariance} Let $X$ and $Y$ be homogeneous graphs. 
Let $Z$ be a graph so that for some $L>0$ :

\begin{itemize}
\item $Z$ contains a disjoint family of subgraphs $\{X_i\}$ and $\{Y_i\}$
whose union contains all the vertices.   
\item Every edge of $Z$ not in one of the subgraphs connects some $X_i$ to
some $Y_j$, and there is exactly one such edge at every vertex of $Z$.
\item For every $i$, there is an $L$ bilipschitz equivalence of $X_i$
 (resp. $Y_i$) and $X$ (resp. $Y$).  
\item The quotient graph obtained from $Z$ by collapsing each of the 
subgraphs to a point is a tree.
\end{itemize}

  There is an $M$, depending only on $L$ and the homogeneity constants
of $X$ and $Y$, so that for any edge, $e$, in $Z$ connecting an $X_i$ and
a $Y_j$ and any choice of base points in $X$ and $Y$, there is an $M$
bilipschitz equivalence of $Z$ to $(X,x_0)*(Y,y_0)$ taking $e$ to the base
edge.
\end{lemma}

\begin{proof}

Note that as $X$ is homogeneous there is a $K$ so that for any
$X_i$, any $x \in X_i$, and any $x' \in X$, there is a $K$ bilipschitz
equivalence $X_i \to X$ which takes $x$ to $x'$.  The same holds for
$Y_j$'s mapping to $Y$.

Call the edge $e$ the base edge of $Z$.  Choose bilipschitz equivalences
of $X_i \to X$  and $Y_j \to Y$, as above.  As $X$ and $Y$ are
homogeneous we may assume that the endpoints of the base edge are $x_0$
and $y_0$.  This  gives a quasi-isometry between the union, $\Sigma_0$, of
$X_i$, $Y_j$, and the base edge to the subgraph of $X*Y$ which was called
$\Gamma_0$ in the earlier construction.   This map restricts to a
bijection of the edges of $(X,x_0)*(Y,y_0)$, not in $X$ or $Y$, with an
endpoint in
$\Gamma_0$ and the edges of $Z$, not in $X_i$ or $Y_j$, with an endpoint
in $\Sigma_0$.

One now follows the construction of $X*Y$.  At every stage we have a
subgraph $\Sigma_n$ of $Z$ with a quasi-isometry to $\Gamma_n$ which
induces a bijection of the incident edges not in the distinguished
subgraphs.  Let$\Sigma_{n+1}$ be the subgraph of $Z$ which contains
$\Sigma_n$, all the $X_i$ and $Y_j$ which are adjacent to it, and the edges
connecting them.  Extend the map to a bilipschitz equivalence of
$\Sigma_{n+1}$ and $\Gamma_{n+1}$, by choosing, for each new $X_i$ and
$Y_j$ a bilipschitz equivalence with the copy of $X$ or $Y$ attached at
the corresponding point of $\Gamma_n$ which takes the point of attachment
to
$x_0$ or $y_0$.  

These bilipschitz equivalences give, in the limit as $n \to \infty$, the
desired bilipschitz equivalence.
  
\end{proof}

  In particular, if $X$ and $Y$ are homogeneous then the bilipschitz type
of $(X,x_0)*(Y,y_0)$ does not depend on the choice of basepoints so we
will usually write simply $X*Y$.  This also implies that $X*Y$ is also
homogeneous.  All the spaces we consider are built out of coset spaces by
passing to bilipschitz equivalent spaces and the free product
construction, and hence are all homogeneous. 

The basepoints within each copy of $X$ or $Y$ in $X*Y$ are, even for
homogeneous $X$ and $Y$, a useful bookkeeping device.  Choosing
basepoints amount to a choice of base edge in $X*Y$; the base points of
each copy of $X$ or $Y$ is determined by being the closest
point in that subgraph to the base edge.

To use Lemma \ref{invariance}, we need constructions of bilipschitz
equivalences.  The next two lemmas are important examples of this, and
clearly demonstrate the utility of the generality of homogeneous spaces
rather than simply coset spaces of groups.

\begin{lemma}\label{plus} Let $X$ be infinite.  Define $X^{+}$ as the graph 
obtained from $X$ by adding a vertex $v$ which is connected by an 
edge to the base point of $X$.  There is a bilipschitz equivalence 
between $X$ and $X^{+}$, hence $X*Y$ and $X^{+}*Y$ are bilipschitz 
equivalent for any $Y$.
\end{lemma}

\begin{proof}

As $X$ is an infinite connected graph of bounded valence there is an 
infinite embedded path in $X$, starting at the base point.  Let 
$x_{0}, x_{1}, \ldots$ be such a path.  Define a map from $X$ to 
$X^{+}$ as follows:

\begin{itemize}
\item Send all points in $X \setminus \{x_{i}\}$ to their images under the 
inclusion of $X$ in $X^{+}$.
\item Send $x_{0}$ to $v$.
\item For $i>0$ send $x_{i}$ to $x_{i-1}$. 
\end{itemize}

It is easy to verify that this is a bilipschitz equivalence.
\end{proof}

Note that this bilipschitz equivalence of $X^+$ and $X$ implies that
if $X$ is homogeneous then so is $X^+$, although it will almost never
have a transitive group of graph automorphisms.

This technique of "sliding from infinity" along a path is used 
repeatedly in the following constructions to produce bilipschitz
equivalences.

Given $X$ and $Y$, define the wedge of $X$ and $Y$ as the space obtained 
from the disjoint union by adding an edge connecting the base points. 
Notice that even if $X$ and $Y$ are homogeneous the wedge, in general,
is not.  Thus the choice of a base point is nontrivial.  We will pick
one of the endpoints of the edge joining the halves.

One important reflection of the fact that free products are tree-like 
is the following:

\begin{lemma}\label{wedge} Let $X$ and $Y$ be infinite homogeneous
space, and let $Z=X*Y$.  There is a bilipschitz equivalence between
$Z$ and the wedge of two copies of $Z$.
\end{lemma}

\begin{proof}
The base edge divides $Z$ into two infinite, connected 
subgraphs.  Thus $Z$ is bilipschitz equivalent to the wedge of these 
two halves.  Each half is almost $X*Y$ - it is covered by disjoint 
families of copies of $X$ and $Y$, connected alternately in a tree of 
spaces.  The way in which the halves differ from $X*Y$ is that in
a single subgraph (isomorphic to $X$ in one half and to $Y$ in the 
other) the base point is not connected to any other subgraph.  

Consider the half in which the deficient subgraph is $X$.  In that copy 
of $X$, choose a path, $p$, from the base point to infinity.  Let
$e_i$, for $i>0$, be the edges connecting $p(i)$ to $y_i$ in copies of
$Y$.  Modify the graph structure by removing the edges $e_i$ and adding
edges $e'_i$ which connect $p(i-1)$ to $y_i$.  This does not change the
bilipschitz type, and, by Lemma \ref{invariance}, the resulting graph is
bilipschitz equivalent to $X*Y=Z$.

After the analogous modification in the other half, the graph is
bilipschitz equivalent to $Z$ wedge $Z$, completing the proof.

\end{proof}

Note that this implies that for $X$ and $Y$ homogeneous, and $Z=X*Y$,
that a wedge of two (and hence any finite number) of copies of $Z$ is
homogeneous.  

As we observe in the introduction, it is not true, even for groups, that 
if $G$ and $G'$ are quasi-isometric then $G*H$ and $G'*H$ are 
quasi-isometric.  Thus one is motivated to ask, as in \cite{G}, when 
quasi-isometric groups are bilipschitz equivalent.  In \cite{P} it is 
shown that all non-abelian free groups are bilipschitz equivalent.  
In \cite{W}, where the general question of when a quasi-isometry is 
at bounded distance of a bilipschitz map is resolved, it is shown that 
any two quasi-isometric non-amenable groups are bilipschitz 
equivalent.  No example of infinite groups which are quasi-isometric 
but not bilipschitz equivalent is known, but \cite{BK} and \cite{McM} 
show there are graphs quasi-isometric to $\Z^{2}$ which are not 
bilipschitz equivalent to $\Z^{2}$.

Consider the special case of Theorem \ref{FreeProducts} where $A$ is a 
subgroup of finite index of $B$.  There is a natural homomorphism 
from $A*C$ to $B*C$, but unless $C$ is trivial, it has infinite index 
image.  There is a subgroup of $B*C$ of finite index built out of $A$ 
and $C$; it is isomorphic to $A*C*\ldots*C$ with $[A:B]$ copies of 
$C$.  Thus, Theorem \ref{FreeProducts} implies that $A*C$ and $A*C*C$ 
are quasi-isometric.  That special case, generalized to arbitrary 
spaces, is one of the key constructions in the proof of the later 
theorems.

\begin{lemma}\label{main} If $X$ and $Y$ are infinite homogeneous spaces 
 then $X*Y$ and $X*(Y*Y)$ are bilipschitz equivalent. 
\end{lemma}

\begin{proof}

We are ready to prove lemma \ref{main}.  By lemma \ref{plus}
$X*Y$ and $X^{+}*Y$ are quasi-isometric. Here we consider
$X^{+}$ as a space with distinguished base point the vertex $x_0$ which
is the base point of the copy of $X$ inside it, and we let $v$ be the
vertex in $X^+ \setminus X$.

  We claim $X^{+}*Y$ is quasi-isometric to $X*(Y*Y)$.  We prove this by 
modifying the graph structure of $X^{+}*Y$ without changing the 
quasi-isometry type, and then verifying that the modified graph satisfies 
the properties which characterize $X*(Y*Y)$ up to quasi-isometry.

We first modify $X^{+}*Y$ as follows: in every copy of $X^{+}$, the 
vertex $v$ is connected to a base point $y_{0}$ of a copy of $Y$ and 
to the base point $x_{0}$ of the copy on $X$ within $X^{+}$.  The 
vertex $x_{0}$ is also connected to a vertex $y_{1}$ in a copy of $Y$.
We add an edge connecting $y_{0}$ to $y_{1}$ and delete the edge 
connecting $v$ to $y_{0}$.  Call this modified graph $G_{0}$.

We have not changed the vertex set, and the identity map on vertices 
is bilipschitz between $X^{+}*Y$ and $G_{0}$.  Further, $G_{0}$ is 
vertex-covered by subgraphs identified with $X^{+}$ and $Y$.  In each 
copy of $X^{+}$ the $v$ vertex is now only connected to $x_{0}$, and not
to any copies of $Y$.  We wish to remove these vertices.  Consider 
the full subgraph of $G_{0}$ which consists of all the vertices except these 
$v$ vertices and their incident edges.  Call this subgraph $G_{1}$.

As $G_{1}$ is connected, and every vertex of $G_{0}$ is in, or 
adjacent to, $G_{1}$, the inclusion of $G_{1}$ into $G_{0}$ is a 
quasi-isometry.  The graph $G_{1}$ consists of copies of $X$ and 
copies of $Y$ and edges connecting them.   

If we delete all edges of $G_1$ that connect a copy of $X$ to a copy
of $Y$, what remains divides into two types of components: copies of $X$
and components that are made up of copies of $Y$'s joined by edges. We
call the latter components of $Y$'s. 

We need to understand the pattern in which the copies of $Y$ are 
connected in a component of $Y$'s.  The edges connecting copies of
$Y$ can be described as follows: for every copy of $X^{+}$ in $X^+*Y$ 
we have, in $G_1$, an edge connecting the two copies of $Y$ that are 
connected in $X^{+}*Y$ to the base point and the "extra" point, $v$. 
In $X^{+} * Y$, every non base point in a copy of $Y$ connects to the 
base point of an $X^+$.  Thus, every non base point in a copy of $Y$ in 
$G_1$ connects to another copy of $Y$.

We now "slide" edges in the components of $Y$ in $G_1$ so that every
vertex connects to another copy of $Y$.  For any component of $Y$ in
$G_1$ whose base point does not connect to another copy of $Y$, choose
a simple path, $p$, to infinity in that $Y$, starting at the base point. 
Modify the edge structure of $G_1$ as follows: for $i>0$ let $e_i$ be the
edge at $p(i)$ connecting to another copy of $Y$, and let $y_i$ be the
other endpoint of $e_i$.  Remove all the $e_i$, and add edges $e'_i$
connecting
$p(i-1)$ to $y_i$. Let $G_2$ be the resulting graph.
 
It is clear that $G_2$ is bilipschitz equivalent to $G_1$.  The graph
$G_2$ now contains subgraphs isomorphic to $X$ and subgraphs bilipschitz
equivalent to $Y*Y$ (in fact, isomorphic to $Y*Y$ aside from
choices of basepoints, which does not change the bilipschitz type, by
Lemma \ref{invariance}), and edges connecting them.  The only difference
between $G_2$ and $X*(Y*Y)$ is  that not every point in each
$Y*Y$ connects to a copy of $X$: the base points of those components
which, in $X^{+}*Y$, connect to the "extra" vertex $v$ in copies of $X^+$
do not connect to copies of $X$ in $G_2$.  

To fix this we again slide edges.  In every copy of $Y$ in $G_2$ whose
base point does not connect to a copy of $X$, choose a $p$ starting at
the base point and running to infinity.  Now let $e_i$, for $i>0$, be the
edges connecting $p(i)$ to $x_i$ in copies of $X$.  Let $G_3$ be the
graph constructed from $G_2$ by removing the $e_i$ and adding $e'_i$
which connects $p(i-1)$ to $x_i$.

 The graph $G_3$ consists of copies of $X$ and of $Y*Y$ connected 
as in Lemma \ref{invariance}, and so is bilipschitz equivalent to
$X*(Y*Y)$.  It is, by construction, bilipschitz equivalent to
$X^{+}*Y$ and therefore to $X*Y$, completing the proof.

\end{proof}

  As discussed above, the lemma is false for some cases of finite $X$ or
$Y$.  It is not difficult, using the above techniques, to determine the
truth in that case, but the answer is somewhat complicated.   For the
cases arising from groups, this is analyzed at the start of the proof of
Theorem \ref{FreeProducts} in the next section.

\section{Building Quasi-Isometries}

In this section we show how to use the constructions of the last 
section to build quasi-isometries between various graphs of groups, 
proving Theorems \ref{FreeProducts}, \ref{Amalgamation}, and 
\ref{General}. 

We start with the proof of Theorem \ref{FreeProducts}.
\begin{proof}

First, if $A$ (and therefore $B$) and $C$ are finite, one knows
that $A*C$ and $B*C$ are virtually free.  Under the assumptions
on cardinality, they are not virtually cyclic, and hence are 
quasi-isometric.

Second, if $A$ and $B$ are finite, but $C$ is infinite, then $A*C$ 
and $B*C$ contain subgroups of finite index isomorphic to 
$C*\ldots*C$ ($|A|$ factors) and $C*\ldots*C$ ($|B|$ factors) 
respectively.  These are quasi-isometric by lemma \ref{main}.

Likewise, if $A$ and $B$ are infinite, but $C$ is finite, then $A*C$ 
and $B*C$ contain subgroups of finite index isomorphic to $A*\ldots*A$ 
and $B*\ldots*B$, both with $|C|$ factors.  Again, by lemma \ref{main} 
$A*\ldots*A$ is quasi-isometric to $A*A$ and $B*\ldots*B$ is quasi-isometric
to $B*B$.

Therefore it suffices to prove the theorem when $A$, $B$, and $C$ infinite. 
We assume now that this is the case. By lemma 
\ref{main}, it suffices to prove that $A*C*C$ and $B*C*C$ are 
quasi-isometric.  So we assume, from now on, that $C$ splits as a 
free product of infinite groups.

Let $f:A \to B$ be a quasi-isometry.  There are 
nets $X$ in $A$ and $Y$ in $B$ so that $f$ induces a bilipschitz 
equivalence $X \to Y$. We can choose such $X$ and $Y$ to include the
base points of $A$ and $B$. Let $r_{1}:A \to X$ and $r_{2}:B \to Y$ be 
projections onto the nets, moving points a uniformly bounded 
distance.  Choose these projections so that only the base point of $A$
maps to the base point of $X$, and likewise only the base point of $B$ 
maps to the base point in $Y$. 

Consider the space $A*C$.  Inside of each copy of $A$, for each $a$, 
there is an edge $e$ which connects $a$ to some $c$ in a copy of 
$C$.  Remove that edge, and replace it by an edge connecting $c$ to 
$r_{1}(a)$.  Since the distance between $a$ and $r_{1}(a)$ is 
uniformly bounded, the new graph is bilipschitz equivalent to $A*C$.  
All the edges leaving each copy of $A$ do so at a point of $X$, hence we 
can replace each copy of $A$ by a copy of $X$ without changing the 
quasi-isometry type.

What we now have is not quite $X*C$, since each point of $X$ connects 
to possibly more than one copy of $C$.  For each $x$ in a copy of $X$, 
pick one of the copies of $C$ connected to $x$, and slide all the 
edges connecting $x$ to other copies of $C$ to connect to the chosen 
copy of $C$ instead of $x$.  Note that if there is more that one copy 
of $C$ connected to $x$, $x$ is not the base point of $X$, and so 
connects to the base point of the copies of $C$.

The resulting graph is a tree of spaces, with copies of $X$ 
connecting to spaces which are made of copies of $C$ attached
to each other by edges joining their base points, in other words, 
wedges of finite number of copies of $C$.  Since $C$ splits 
as a free product of infinite groups, Lemma \ref{wedge} shows each of 
these wedges is bilipschitz equivalent to $C$. Thus $A*C$ 
is quasi-isometric to $X*C$ by Lemma \ref{invariance}.

By the same construction $B*C$ is quasi-isometric to $Y*C$.  Since $X$
and $Y$ are bilipschitz equivalent, $X*C$ and $Y*C$ are bilipschitz
equivalent, which completes the proof.  We note that the proof goes 
through unchanged for homogeneous spaces rather than groups.

\end{proof}

Stallings' Ends theorem says that any group with infinitely many ends 
splits non-trivially over a finite group, thus we want to extend 
Theorem \ref{FreeProducts} to cover such splittings.  This is the 
content of Theorem \ref{Amalgamation}, which we now prove.

\begin{proof}

To start, consider a free product with amalgamation $A*_{F}B$, with 
$F$ a finite normal subgroup of both $A$ and $B$.  In this case we 
have $F$ as a normal subgroup of $A*_{F}B$, with quotient 
$(A/F)*(B/F)$.  Thus $A*_{F}B$ is quasi-isometric to $(A/F)*(B/F)$.  
Since $A/F$ is quasi-isometric to $A$ and $B/F$ is quasi-isometric to 
$B$, Theorem \ref{FreeProducts} proves $A*B$ and $A*_{F}B$ are 
quasi-isometric.

When $F$ is not normal in $A$ or $B$, we make the same argument, but 
now $A/F$ and $B/F$ are spaces rather than groups. In the tree 
of spaces modeling $A*_{F}B$ (see \cite{SW}), one has copies of $A$ 
and $B$, but rather than single edges connecting copies of $A$ to 
$B$, one has an $F$ coset in $A$ joined to an $F$ coset in $B$ by 
$|F|$ edges.  

Choose $X$ in $A$ so that $X$ contains one point of each coset $aF$. 
Give $X$ the structure of a graph by joining $x$ and $x'$ by an edge if
there is an edge joining points in the corresponding cosets.  Likewise,
choose $Y$ in $B$ which intersects every $F$ coset in one point. Replacing
each $A$ by $X$ and $B$ by $Y$, with an edge joining
$x \in X$  and $y \in Y$ if and only if the corresponding cosets are
connected in $A*_F B$, gives a quasi-isometry between $A*_{F}B$ and
$X*Y$.  As $X$ is quasi-isometric  to $A$ and $Y$ is quasi-isometric to
$B$, it follows exactly as in the proof of Theorem \ref{FreeProducts}
that $A*B$ is quasi-isometric to $X*Y$.  Thus $A *_F B$ and $A*B$ are
quasi-isometric.

 The situation for $A*_{F}$ is similar.  The model for
$A*_{F}$ is a tree of spaces, each of which is copy of $A$.  The edges
between copies of $A$  are directed, attaching a coset of the first
embedding of $F$ in $A$ in one  copy of $A$ to a coset of the second
embedding of $F$ in $A$ in another copy of $A$.  Every point is connected
to two other copies of $A$, once as an  initial vertex and once as a
terminal vertex.

To carry out the same argument, one needs to find a subset $X$ of $A$ 
which is simultaneously a set of coset representatives for both 
embeddings of $F$ in $A$.  That such an $X$ exists is a standard 
application of Hall's Marriage Lemma (\cite{Hall}).  Given such an $X$,
the argument  above.  
\end{proof}

Having proven Theorems \ref{FreeProducts} and \ref{Amalgamation}, it 
is straightforward to deal with an arbitrary graph of groups with finite 
edge groups, as such a graph of groups is simply iterated free products 
with finite amalgamation and HNN extensions over finite subgroups (\cite{S}).

In order to prove Theorem \ref{General}, we prove a slightly 
different result which clearly implies it, but is somewhat more
awkward.

\begin{theorem} Let $G$ be a graph of groups with finite edge groups.  
Let $S$ be the set of quasi-isometry types of vertex spaces, without 
repetition.  Let $G'$ be the free product of a finite set of groups 
with quasi-isometry types representing every type of $S$ exactly 
once.  Let $F$ and $F'$ be any (possibly trivial or cyclic) free groups
so that $G*F$ and $G'*F'$ have the same number of ends, then $G*F$ 
and $G'*F'$ are quasi-isometric.
\end{theorem}

\begin{proof}
We prove this by induction on the number of edges in the graph.  If 
there are no edges the result is essentially a tautology, given that 
one knows the number of ends is a quasi-isometry invariant which
classifies free groups up to quasi-isometry.  

As a graph of groups with $n+1$ edges can be built out of graphs with 
fewer edges, either by free product or HNN extension, the result 
follows from the following, which is immediate from Theorems 
\ref{FreeProducts} and \ref{Amalgamation}.

\begin{theorem}
If $A$ is quasi-isometric to $A'$ and $B$ is quasi-isometric to 
$B'$, and $F$ a finite proper subgroup of $A$ and $B$, then $A*_{F}B$ is 
quasi-isometric to $A'*B'$ unless both products are virtually free. Likewise 
$A*_{F}$ is quasi-isometric to $A'*\Z$ unless both are virtually free.
\end{theorem}  
\end{proof}

\section{Obstructing Quasi-Isometries}

In this section we prove Theorem \ref{Accessible}, a partial converse to 
the earlier theorems.  In view of Theorem \ref{General}, and the fact 
that accessibility is a quasi-isometry invariant (\cite{TW}), this 
comes down to:

\begin{theorem}\label{Factor}
Let $G_{1}$ and $G_{2}$ be quasi-isometric groups, both of 
which are fundamental groups of terminal graphs of groups.   If $H$ 
is a one-ended vertex group of $G_{1}$ then there is a  one-ended vertex
group of $G_{2}$ which is quasi-isometric to $H$.
\end{theorem}
\begin{proof}

\begin{lemma} Let $G$ be the fundamental group of a terminal graph of
groups, and let $H$ be a one-ended group.  For any $(A,B)$ there is an
$R$ so that for any $(A,B)$ quasi-isometric embedding $f$ of $H \to G$
there is a (necessarily unique) vertex space $X$ in $G$ with $H$ contained
in the $R$ neighborhood of $X$.
\end{lemma}

\begin{proof} 

Let $D$ be the maximal diameter of an edge space of $G$.  For any such
edge space $E$, the pre-image $f^{-1}(E)$ in $H$ has diameter at most
$A(D+B)$.   By the definition of one-endedness, there is some $D'$ so that
all but one of the components of the complement of any set $S$ of
diameter at most $A(D+B)$ lie entirely within $D'$ of $S$.

Thus, for any edge space $E$, $f(H)$ is contained in the
$AD'+B+1$ neighborhood of one side of $E$.  If $f(H)$ is not contained
within twice this distance of any vertex, then we can orient every
edge to point towards the half containing $f(H)$, and there is at least
one edge pointing away from every vertex.  Thus there are unbounded
oriented rays.  On the other hand, if $v$ is a vertex space which has
nontrivial intersection with $f(H)$, then every edge more than $AD'+B+1$
from $v$ must be oriented towards $v$, which contradicts the existence of
unbounded oriented rays.  Thus the hypothesis that $f(H)$ is not
contained in a $2(AD'+B+1)$ neighborhood must be false.

\end{proof}

Theoem \ref{Factor} now follows easily.  The lemma shows $f(H)$ must be
contained in a neighborhood of a vertex space $K$ of $G_2$.   Applying
the lemma to the inverse quasi-isometry, $f'$, gives $f'(K)$ contained in
a neighborhood of some $H'$ in $G_1$.  Thus $f'f(H)$ is contained in a
neighborhood of $H'$.  Since $f'f$ is within bounded distance of the
identity map, this implies $H$ is contained in a neighborhood of $H'$
which implies $H=H'$.  Further, the fact that a neighborhood of $f'f(H)$ 
contains $H$ implies that a neighborhood of $f(H)$ contains $K$, so $f$
restricts to give a quasi-isometry between $H$ and $K$.

\end{proof} 

Theorem \ref{Accessible} reduces the large scale geometry of 
accessible groups to the large scale geometry of one-ended groups.  
It would be very interesting to understand the geometry of 
non-accessible groups.

\noindent
Panos Papasoglu\\
Departement de Math\'ematiques, Universit\'e de Paris XI (Paris-Sud)\\
91405 Orsay, FRANCE\\
E-mail: panos@math.u-psud.fr\\

\noindent
Kevin Whyte\\
Dept. of Mathematics\\
University of Chicago\\
Chicago, Il 60637\\
E-mail: kwhyte@math.uchicago.edu\\


\begin{thebibliography}{}

\bibitem[BK]{BK}
D. Burago and B. Kleiner, Separated nets in Euclidean Space, preprint

\bibitem[D1]{D1}
M.J.Dunwoody, The accessibility of finitely presented groups, Invent. 
Math., p. 449-457, 1985

\bibitem[D2]{D2} 
M.J.Dunwoody, An inaccessible group, The Proceedings of Geometric 
Group Theory 1991, G.A.Niblo, M. Roller (Eds)
LMS Lecture Notes Series 181, Cambridge University Press,p.75-78, 1993

\bibitem[G]{G}
M.Gromov, Asymptotic invariants of infinite groups, in 'Geometric group theory', 
(G.Niblo, M.Roller, Eds.), LMS Lecture Notes, vol. 182, Cambridge Univ. 
Press, 1993

\bibitem[GW]{Hall}
J. Graver and M. Watkins, Combinatorics with Emphasis on the Theory of
Graphs, Springer-Verlag, 1977

\bibitem[McM]{McM}
C. McMullen, Lipschitz maps and nets in Euclidean space, preprint

\bibitem[P]{P} 
P.Papasoglu, Homogeneous Trees are Bilipschitz Equivalent, 
Geometriae Dedicata, vol. 54, p. 301-306, 1995

\bibitem[S]{S}
J-P.Serre, Trees, Springer Verlag, N.Y., 1980 

\bibitem[St1]{St1}
J.R.Stallings, On torsion-free groups with infinitely many ends, 
Ann. of Math. 88, p. 312-334, 1968 

\bibitem[St2]{St2}
J.R.Stallings, Group theory and three dimensional manifolds, Yale 
Mathematical Monographs 4, Yale University Press, New Haven, 1971

\bibitem[SW]{SW}
P. Scott and T. Wall,Topological methods in group theory, London Math. Soc. Lecture Notes Ser. 36, pp.137-203. Cambridge University Press, 1979.

\bibitem[TW]{TW}
C.Thomassen and W.Woess, Vertex transitive graphs and accessibility,
J. Comb. Theory, Ser.B 58, No.2, p.248-268, 1993

\bibitem[W]{W}
K. Whyte, Amenability, Bilipschitz Equivalence, and the Von Neumann 
Conjecture, Duke Mathematics Journal, vol. 99, p.93-112, 1999

\end{thebibliography}
\end{document}